\input amstex
\documentstyle{amsppt}
\voffset=-3pc
\def\b1{\bold 1}
\def\bC{\Bbb C}

\def\bM{\Bbb M}
\def\ck{\Cal K}

\def\ob{\overline \beta}

\def\ep{\epsilon}

\def\ssa{S_{\text{sa}}}
\def\asa{A_{\text{sa}}}
\def\bsa{B_{\text{sa}}}

\def\tA{\widetilde A}
\def\th{\widetilde h}

\def\Asa{\widetilde A_{\text{sa}}}

\magnification=\magstep1
\parskip=6pt
\NoBlackBoxes
\topmatter
\title Some Directed Subsets of C*--algebras and Semicontinuity Theory
\endtitle  
\author Lawrence G.~Brown
\endauthor

\abstract{The main result concerns a $\sigma-$unital $C^*$--algebra $A$, a strongly lower semicontinuous element $h$ of $A^{**}$, the enveloping von Neumann algebra, and the set of self--adjoint elements $a$ of $A$ such that $a\le h-\delta\b1$ for some $\delta >0$, where {\bf 1} is the identity of $A^{**}$. The theorem is that this set is directed upward. It follows that if this set is non-empty, then $h$ is the limit of an increasing net of self--adjoint elements of $A$. A complement to the main result, which may be new even if $h=\b1$, is that if $a$ and $b$ are self--adjoint in $A$, $a\le h$, and $b\le h-\delta\b1$ for $\delta>0$, then there is a self--adjoint $c$ in $A$ such that $c\le h, a\le c$, and $b\le c$.}
\endabstract
\endtopmatter

Let $A$ be a $C^*$--algebra and denote by $A^{**}$ its bidual, the enveloping von Neumann algebra. In [AP] C. Akemann and G. Pedersen discussed several classes of self--adjoint elements of $A^{**}$ as analogues of semicontinuous functions on topological spaces.

For a subset $S$ of $A^{**}$, $\ssa$ denotes $\{h\in S: h^*=h\}$, $S_+$ denotes $\{h\in S; h\ge 0\}$, and $S^-$ denotes the norm closure of $S$. And for $S\subset A^{**}_{\text {sa}}$, $S^{\text m}$ denotes the set of ($\sigma$-weak) limits of bounded increasing nets in $S$, $S_{\text m}$ denotes the set of limits of bounded decreasing nets in $S$, and $S^\sigma$ denotes the set of limits of bounded increasing sequences in $S$. Also $\tA$ denotes $A+\bC${\bf 1}, where {\bf 1} is the identity of $A^{**}$, $LM(A)=\{T\in A^{**} : TA\subset A\}, QM(A)=\{T\in A^{**} : ATA\subset A\}$, $\ck$ is the $C^*$--algebra of compact operators on a separable infinite dimensional Hilbert space, and $\bM_2$ is the $C^*$--algebra of $2\times 2$ matrices. Finally, $S(A)$ denotes the state space of $A$, and states of $A$ are also considered as (normal) states of $A^{**}$.

Three of the classes considered in [AP] are $(A^{\text m}_{\text {sa}})^-$, whose elements are called strongly lsc, $(\Asa)^{\text m}$, whose elements are called middle lsc, and $((\Asa)^{\text m})^-$, whose elements are called     weakly lsc. Also $h$ is called usc in any sense if $-h$ is lsc in that sense. For example, $((\Asa)_{\text m})^-$ is the set of weakly usc elements. Akemann and Pedersen showed in [AP] that these three kinds of semicontinuity can all be different. But they also considered the class $\asa^{\text m}$, and left open the question whether $A_{\text {sa}}^{\text m}$ always equals $(\asa^{\text m})^-$.

In [B, Corollary 3.25] I showed that $\asa^{\text m}$ does equal $(\asa^{\text m})^-$ if $A$ is separable. Whether this is so for general $A$ is still open, and I have no conjecture concerning this. I think the main obstacle to proving the affirmative answer is that it is difficult to construct increasing nets in non-commutative $C^*$--algebras. I know two ways of doing this. In cases where sequences suffice one can use recursive constructions. Or one can find a subset of $A_{\text {sa}}$ which is a directed set in the natural ordering. The first method was used in the proof of [B, Corollary 3.25], and the second method has been used many times. For example, it is used in the construction of approximate identities of $C^*$--algebras, and it is used several times in [AP]. It seems to me that the best hope for proving $A^{\text m}_{\text {sa}}=(\asa^{\text m})^-$ when $A$ is non-separable is to prove that a suitable subset of $A_{\text {sa}}$ is directed upward.

For $h$ in $(\asa^{\text m})^-$ consider the following sets:

Let $\Cal A_1=\{a\in A_{\text {sa}}:\quad a\le h\}$

Let $\Cal A_2=\{a\in A_{\text {sa}}:\quad \varphi(a)<\varphi(h), \forall\varphi\in S(A)\}.$

If $h\ge 0$, let $\Cal A_3=\{a\in A_+: a\le\theta h$ for some $\theta$ in $(0,1)\}$.

If $\varphi(h)>0, \forall\varphi\in S(A)$, let
$\Cal A_4=\{a\in A_+:\exists \theta$ in $(0,1)$ such that $\varphi(a)<\varphi(\theta h),\forall\varphi\in S(A)\}$.

Let $\Cal A_5=\{a\in A_{\text {sa}}: a\le h-\delta\b1$ for some $\delta>0\}$.\newline
Examples will be given at the end of this note to show that, even if $A$ is separable, none of $\Cal A_1, \Cal A_2, \Cal A_3, \Cal A_4$ need be directed upward. Counterexamples for $\Cal A_1$ were already given in [B, 3.23 and 5.6]. The main result of this note is that $\Cal A_5$ is directed upward if $A$ is $\sigma$-unital. However, this does not imply that $\asa^{\text m}=(\asa^{\text m})^-$ for  $A$ $\sigma$-unital, since there are some choices of $h$ for which $\Cal A_5$ is empty.

\proclaim{Proposition 1}Let $A$ be a $C^*$--algebra and $a,b\in A_+$ such that $\|a\|\le 1$ and $\|b\|< 1$. Then there is $c$ in $A$ such that $a\le c, b\le c$, and $\|c\|\le 1$.
\endproclaim

\demo{Proof}Let $c=b+(\b1-b)^{\frac 12}[(\b1-b)^{-\frac 12}(a-b)(\b1-b)^{-\frac 12}]_+(\b1-b)^{\frac 12}$. Since \newline
$(\b1-b)^{-\frac 12}(a-b)(\b1-b)^{-\frac 12}\le (\b1-b)^{-\frac 12}(\b1-b)(\b1-b)^{-\frac 12}=\b1$, then $0\le t\le \b1$ if $t=[(\b1-b)^{-\frac 12}(a-b)(\b1-b)^{-\frac 12}]_+$. Therefore $0\le (\b1-b)^{\frac 12}t(\b1-b)^{\frac 12}\le (\b1-b)$, and hence $b\le c\le\b1$. Also since $t\ge (\b1-b)^{-\frac 12}(a-b)(\b1-b)^{-\frac 12}$, then $c\ge b+(a-b)\ge a$. 
\enddemo

\proclaim{Lemma 2} Let $A$ be a separable stable $C^*$--algebra, $h\in (A^{\text m}_{\text {sa}})^-$, and let $\Cal A=\{a\in A_{\text {sa}}: h-a\ge\delta\b1$ for some $\delta>0\}$. Then $a$ is directed upward. Also if $a\in A_{\text {sa}}, a\le h$, and $b\in \Cal A$, then there is $c$ in $A_{\text {sa}}$ such that $c\le h, a\le c$, and $b\le c$.
\endproclaim

\demo{Proof}Assume $\Cal A\ne\emptyset$. Since the hypothesis and conclusion are unaffected if $h$ is replaced by $h-a$, $a\in A_{\text {sa}}$, we may assume $h\ge \ep\b1$ for some $\ep >0$. By [B, Theorem 4.4(a)] $h=TT^*$ for some $T$ in $LM(A)$. Therefore $T^*T\in QM(A)$ and $\sigma(T^*T)$ omits the interval $(0,\ep)$. Let $p$ denote the range projection of $T^*T$. Then $p$ is closed (cf. [B, Proposition 2.44(b)]). Let $T'=(T^*T)^{-1}T^*$, where the inverse is computed in $pA^{**}p$. Then $TT'=\b1$ and $T'T=p$.

We claim that $TA=A$. To see this, note first that $T=Tp$ and hence $TA=(Tp)(pA)$. Since $pA$ is norm closed by [E] and $(Tp)^*(Tp)\ge\ep p$, it follows that $TA$ is norm closed in $A$. Therefore it is sufficient to show that $TA$ is $\sigma$-weakly dense in $A^{**}$. But the $\sigma$-weak closure of $TA$ includes $TA^{**}$, and $TA^{**}\supset TT'A^{**}=A^{**}$.

Now it is clear that $TAT^*=(TA)(AT^*)=(TA)(TA)^*=A\cdot A=A$. If $TaT^*\in \Cal A$, then since $TT^*\ge TaT^*+\delta\b1$, we have $T'(TT^*){T'}^{*}\ge(T'T)a(T'T)^*+\delta T'{T'}^{*}=pap+\delta(T^*T)^{-1}$. Since $T'(TT^*){T'}^{*}=p$ and $(T^*T)^{-1}\ge \eta p$ for some $\eta >0$, this implies $p\ge pap+\delta'p$ for some $\delta'>0$. Conversely, if $p\ge pap+\delta'p$ for $\delta'>0$, then 
$$\align
h=TT^*=TpT^*&\ge T(pap)T^*+\delta'(TT^*)=TaT^*+\delta'h\\
 &\ge TaT^*+\delta'\epsilon\b1.\endalign$$
Therefore $TaT^*\in \Cal A$. By [B, Corollary 3.4] $\{x\in pA_{\text {sa}}p: x\le (1-\delta')p$ for some $\delta'>0\}=\{pyp: y\in A_{\text {sa}}\,\, \text{and}\,\, y\le(1-\delta')\b1$ for some $\delta'>0\}$. It is known that the last set is directed upward (cf. [P2, Theorem 1.4.2]), so it follows that $\Cal A$ is directed upward. The last sentence is proved similarly using Proposition 1.
\enddemo

\proclaim{Lemma 3} Let $A$ be a separable $C^*$--algebra, $h\in (A^{\text m}_{\text {sa}})^-$, and let $\Cal A=\{a\in A_{\text {sa}}: h-a\ge\delta\b1$ for some $\delta>0\}$. Then $\Cal A$ is directed upward. Also if $a\in A_{\text {sa}}, a\le h$, and $b\in \Cal A$, then there is $c$ in $A_{\text {sa}}$ such that $c\le h, a\le c$, and $b\le c$.
\endproclaim

\demo{Proof}As in the proof of Lemma 2, we may assume $h\ge \ep \b1$ for $\ep>0$. Consider $B=A\otimes \ck$ and $\th=h\otimes\b1$ in $B^{**}$. Obviously $\th\in (B_{\text {sa}}^{\text m})^-$. Let $\widetilde{\Cal A}=\{b\in B_{\text {sa}}: \th-b\ge \delta\b1$ for some $\delta>0\}$, and let $e$ be a rank one projection in $\ck$. If $a_1, a_2\in \Cal A$, then $a_1\otimes e, a_2\otimes e\in \widetilde{\Cal A}$. By Lemma 2, $\exists b\in B_{\text {sa}}$ and $\delta>0$ such that $a_i\otimes e\le b\le\th-\delta\b1$ for $i=1,2$. Therefore $a_i\otimes e\le (\b1\otimes e) b(\b1\otimes e)\le (\b1\otimes e)\th (\b1\otimes e)-\delta (\b1\otimes e)$. Under the canonical isomorphism of $A$ with $(\b1\otimes e)B(\b1\otimes e)$, this is the desired conclusion. The last sentence is proved similarly. 
\enddemo

The following is implicit in the proof of [B, Theorem 3.24], but we provide a proof here.

\proclaim{Lemma 4} If $A$ is a  $C^*$--algebra, $h\in (A^{\text m}_{\text {sa}})^-$, and if $\{a_n\}$ is a countable subset of $A_{\text {sa}}$ such that $a_n\le h, \forall n$,  then there is $h_1\in A^\sigma_{\text {sa}}$ such that $h_1\le h$ and $a_n\le h_1, \forall n$.
\endproclaim

\demo{Proof} We construct recursively an increasing sequence $(b_n)$ in $A_{\text {sa}}$ such that $a_k\le b_n+\frac 1n\b1$ for $1\le k\le n$ and $b_n\le h,\forall n$. Then let $h_1=\lim b_n$. Let $b_1=a_1$. If $n>1$ and $b_1,\dots, b_{n-1}$ have been constructed, then by [B, Corollary 3.17] there is a function $f$ such that $\lim\limits_{\ep\to 0^+} f(\ep)=0$ and the following is true: If $\ep>0, x\in A_{\text {sa}}$, and $b_{n-1}-\ep\b1\le x\le h+\ep\b1$, then there is $y$ in $A_{\text {sa}}$ such that $b_{n-1}\le y\le h$ and $\|y-x\|<f(\ep)$. Choose $\ep$ in $(0,{1/2n})$ such that $f(\ep)< 1/2n$. By [B, Theorem 3.24(b)] there is $x$ in $A_{\text {sa}}$ such that $a_k-\ep\b1\le x\le h$ for $i=1,\dots, k$ and $b_{n-1} -\ep\b1\le x$. Choose $b_n$ to be the $y$ indicated above, then $a_k\le x+{1\over 2n}\b1\le y+{1\over n}\b1$ for $k=1,\dots, n$. 
\enddemo

\proclaim{Theorem 5}Let $A$ be a $\sigma$-unital $C^*$--algebra, $h\in (A^{\text m}_{\text {sa}})^-$, and let $\Cal A=\{a\in A_{\text {sa}}: h-a\ge \delta\b1$ for some $\delta >0\}$. Then $\Cal A$ is directed upward. Also if $a\in A_{\text {sa}}, a\le h$, and $b\in \Cal A$, then there is $c$ in $A_{\text {sa}}$ such that $c\le h, a\le c$, and $b\le c$. 
\endproclaim

\demo{Proof}As in the proof of Lemma 2, we may assume $h\ge\ep\b1$ for some $\ep>0$. Consider $a_1, a_2\in A_{\text {sa}}$ such that $a_i\le h-\delta_i\b1$ for $i=1,2$ and $\delta_i>0$. Let $(e_n)$ be a sequential approximate identity of $A$. Then $\{a_i+\delta_ie_n : i=1,2, n=1,2,\dots\}$ is a countable subset of $\{a\in A_{\text {sa}}: a\le h\}$. Let $h_1$ in $A^\sigma_{\text {sa}}$ be as in Lemma 4. If $h_1=\lim b_m$, $b_m\in A_{\text {sa}}$, let $A_1$ be the separable $C^*$--algebra generated by $a_1, a_2$, the $e_n$'s, and the $b_m$'s. Thus $A^{**}_1$ and $A^{**}$ have the same identity element. Since $a_i+\delta_ie_n\le h_1, n=1,2,\dots$, then $a_i+\delta_i\b1\le h_1$. Therefore by Lemma 3, there is $c$ in $(A_1)_{\text {sa}}\subset A_{\text {sa}}$ such that $a_i\le c$ for $i=1,2$ and $c\le h_1-\delta\b1\le h-\delta\b1$ for some $\delta>0$. The last sentence is proved similarly.
\enddemo

\proclaim{Corollary 6}Let $A$ be a $\sigma$-unital $C^*$--algebra and $h$ an element of $(A^{\text m}_{\text {sa}})^-$ such that $h\ge a+\delta\b1$ for some $a$ in $A_{\text {sa}}$ and some $\delta>0$. Then $h\in A^{\text m}_{\text {sa}}$
\endproclaim

\demo{Proof} The theorem produces an increasing net in $A_{\text {sa}}$ whose limit is the supremum of $\Cal A$ in $A^{**}_{\text {sa}}$. To show that this supremum is $h$, we may assume $h\ge \ep\b1$ for some $\ep>0$, as in the proof of Lemma 2. Then by [AP, Proposition 3.5], $h^{-1}\in((\Asa)_{\text m})^-$.  Next [AP, Proposition 3.1] shows that there is a net $(h_i)$ in $\Asa$ such that $h_i\ge h^{-1}, \forall i$, and $h_i\to h^{-1}\,  \sigma$-weakly. Finally [P1, Lemma 3.2] implies that $h_i^{-1}\to h$ $\sigma$-strongly. (So far we have just used the argument for [AP, Theorem 3.3] in a slightly different context.) Now for each $i$, $h_i^{-1}=\delta_i\b1+a_i$ for some $\delta_i>0$ and $a_i\in A_{\text sa}$. If $e\in A_+$ and $\|e\|<1$, then $\delta_i e+a_i\in \Cal A$, and the supremum of these elements, for fixed $i$, is $h_i^{-1}$. So the supremum of $\Cal A$ is indeed equal to $h$, whenever $\Cal A\ne\emptyset$.
\enddemo

\example{Examples and Remarks}(i) There are a separable unital $C^*$--algebra $A$, a positive $h\in(A^{\text m}_{\text {sa}})^-$, $a_+, a_-\in A_+$ and $\theta\in (0, 1)$ such that $a_i\le\theta h$ and there does not exist $b$ in $A_{\text {sa}}$ with $a_{\pm}\le b\le h$. This shows that neither $\Cal A_1$ nor $\Cal A_3$ is directed upward.

Let $A$ be the algebra of convergent sequences in $\bM_2$. Then $A^{**}$ can be identified with the set of bounded collections, $h=\{h_n\}_{1\le n\le\infty}$ with each $h_n\in\bM_2$. Here $h\in A$ if and only if $h_n\to h_\infty$. Define $h$ in $(A^{\text m}_{\text {sa}})^-$ by $h_n=\pmatrix 1 & 0\\ 0 &\frac 1n\endpmatrix$ for $n<\infty$ and $h_\infty=\pmatrix \frac 12 & 0\\ 0 &0\endpmatrix$. Let $(a_\pm)_n=\pmatrix \theta/2 & \pm n^{-\frac 12} \theta/2\\ \pm n^{-\frac 12} \theta/2 & n^{-1} \theta/2\endpmatrix=\theta h_n^{1\over 2}\pmatrix 1/2 & \pm 1/2\\ \pm 1/2 & 1/2\endpmatrix h_n^{1\over 2}$ for $n<\infty$ and $(a_\pm)_\infty=\pmatrix \theta/2 & 0\\ 0 &0\endpmatrix$, where $\theta$ is to be determined. If $a_\pm\le b\le h$, then $b_n=h_n^{1\over 2} \pmatrix \alpha_n & \beta_n\\ \ob_n &\gamma_n\endpmatrix h_n^{1\over 2}$ where
$$\theta \pmatrix 1/2 & \pm 1/2\\ \pm 1/2 &1/2\endpmatrix \le \pmatrix \alpha_n & \beta_n\\ \ob_n &\gamma_n\endpmatrix \le 1.$$
Note that the larger of $|\beta_n\pm \theta/2|$ is at least $\theta/2$. Therefore $(\theta/2)^2\le(\alpha_n-\theta/2)(\gamma_n-\theta/2)$. Since $\gamma_n\le 1$, then $\alpha_n-\theta/2\ge {\theta^2/4\over 1-\theta/2}$. For suitable $\theta$, this implies $\alpha_n\ge 3/4, \forall n$. Since $\lim \alpha_n\le{1\over 2}$, this is a contradiction.

(ii) It is impossible to have a separable unital counterexample for $\Cal A_2$ or $\Cal A_4$, but a slight variation of the above gives a separable non-unital counterexample. Keeping the above notation, let $B$ be the hereditary  $C^*$--subalgebra of $A$ consisting of sequences with limit of the form $\pmatrix * & 0\\ 0 &0\endpmatrix$. Then $a_\pm\in B$ and $h\in B^{**}\subset A^{**}$. Also $\varphi(h)>0, \forall \varphi\in S(B)$. If $\theta <\theta'<1$, then $\varphi(\theta'h-a_\pm)>0, \forall \varphi\in S(B)$. Therefore neither $\Cal A_2$ nor $\Cal A_4$ is directed upward.

(iii) For an arbitrary $C^*$--algebra $A$, if $h\in(A_{\text {sa}}^{\text m})^-$ and $h\ge 0$, then $q$, the range projection of $h$, is open (cf. [B, Proposition 2.44(a)]). Let $B$ be the hereditary $C^*$--subalgebra supported by $q$. Then $h\in B^{**}\subset A^{**}$, $\varphi(h)>0, \forall \varphi\in S(B)$, and $h\in (\bsa^{\text m})^-$ ([B, Proposition 2.14]). Clearly if $a\in A_+$ and $a\le h$, then $a\in B$.

If $A$ is separable, then it follows from [B, Theorem 3.24(a)] that $h\ge e$ for some strictly positive element $e$ of $B$. Therefore if $a\in \Cal A_4$, then $\theta h-a\ge f$ for some strictly positive element $f$ of $B$. 

There are some variants on $\Cal A_2, \Cal A_3, \Cal A_4$ that I could have included on the list near the beginning of this paper, and this remark shows why I didn't bother.

(iv) The arguments for Lemmas 2 and 3 apply more generally and can prove the following:\newline

Let $A$ be a separable $C^*$--algebra and $h$ a positive strongly lsc element of $A^{**}$ whose spectrum omits $(0,\ep)$ for some $\ep>0$. Then for $B$ as in the previous remark, there are a closed projection $p$ in $(A\otimes\ck)^{**}$ and a complete order isomorphism $\varphi : B \to p(A\otimes\ck)p$ such that $\varphi^{**}(h)=p$. Hence $\varphi^{**}$ is the natural extension of $\varphi$ to the biduals, and in particular for $b$ in $
B_{\text {sa}}$, $b\le h$ if and only if $\varphi(b)\le p$.
\endexample

\bye